
\input amstex
\documentstyle{amsppt}
\pagewidth{125mm}
\pageheight{200mm}

\topmatter
\title
Modified Novikov--Veselov equation and \\
differential geometry of surfaces
\endtitle
\author
Iskander A. TAIMANOV
\endauthor
\address
Institute of Mathematics, 630090 Novosibirsk, Russia
\endaddress
\email
taimanov\@math.nsk.su
\endemail
\endtopmatter
\leftheadtext{}
\rightheadtext{}
\document
\NoBlackBoxes

\head
1. Introduction
\endhead

In the present paper we  consider global soliton deformations of surfaces
immersed in the three-dimensional Euclidean
space.

Local deformation of surfaces represented via the generalized
Weierstrass formulas (3.3 - 3.4) were introduced by Konopelchenko
(\cite{Kon}) by using of the modified Novikov-Veselov
equation (2.11) which in turn was introduced by Bogdanov {\cite{Bg1}).
This equation is a modification of the Novikov--Veselov
equation in the same sense as the modified Korteweg--de Vries equation is
a modification of the Korteweg--de Vries equation.
Saying about the geometric meaning of the Novikov--Veselov equation
we notice that
it has important applications to
the theory of algebraic curves (\cite{T,Sh}).

Here we will discuss applications of the modified Novikov--Veselov
equation to differential geometry of surfaces.
Investigation of global properties of
mNV deformations of surface in the case of tori of revolution was started in
\cite{KT} where it was particularly shown that tori of revolution
are preserved by these deformations. But at that time the relation of
these deformations to conformal geometry was not understood.

In the present paper we consider  global deformations of
surfaces of general type and it's relation to the theory of the
Willmore functional which is defined as integral of squared mean
curvature (5.1).

Locally any regular surface is represented via the generalized
Weierstrass formulas and moreover any analytic surface is globally
represented in this manner
(see Propositions 1 and 2). Thus at least for analytic surfaces the
$mNV$--deformation is correctly defined.

We show that

\vskip3pt

{\it $mNV$-deformation transforms tori into tori and
preserves their conformal structure and value of the Willmore
functional (Theorems 1 and 2).}

\vskip3pt

We also consider the following conjecture

\vskip3pt

\proclaim{Conjecture}
Non-stationary, with respect to $mNV$-deformation, torus
can not be a local minimum of the Willmore functional
\endproclaim

\vskip3pt

\noindent
and discuss it's relation to the famous Willmore conjecture
(see chapter 5).
We thank M.V. Babich, B.G. Konopelchenko, and S.V. Manakov for helpful
conversations.

\head
2. Modified Novikov-Veselov equation
\endhead

\subhead
2.1. Novikov-Veselov equation
\endsubhead

"The Novikov-Veselov equation (NV)
$$
U_t = \partial^3 U + \bar \partial ^3 U + \partial (V U)
+ \bar \partial ( \bar V U),
\eqno{(2.1)}
$$
$$
\bar \partial V = 3 \partial U
$$
was introduced by Novikov and Veselov
in \cite{VN1} within frames of development of the theory of
two-dimensional potential Schr\"odinger operators
which are finite-zone on one level of energy  (\cite{DKN,VN2}).

‡Here functions $U$ and $V$
are defined on the complex plane and  partial derivatives
$\partial , \bar \partial $ are defined by usual formulas
$$
\partial = \frac{1}{2} (\frac{\partial }{\partial x} - i
\frac{\partial }{\partial y}),\ \ \
\bar \partial = \frac{1}{2} (\frac{\partial }{\partial x}
+ i \frac{\partial }{\partial y}),\ \ \
z = x+ iy \in {\bold C}.
$$

ŽThis equation is a natural two-dimensional generalization of the
famous Korteweg-de Vries equation (KDV)
$$
U_t = \frac{1}{4}U_{xxx} + \frac{3}{2} U U_x,
\eqno{(2.2)}
$$
to which the NV equation reduces, after suitable renormalization
of variable $x$,
in the case when the function
$U(z,\bar z)$ does not depend on variable $y$.
Investigation of the KDV equation and discovery of
it's prominent properties
were starting points of intensive development of the soliton theory.
Though the theory of this equation is
well-known due to a lot of monographs
(see, for instance \cite{N}) we will dwell on
some facts which are substantial
for our exposition.

The KDV equation has a representation as a condition of commutativity
$$
[L^{KDV}, \frac{\partial }{\partial t} - A] = 0,
\eqno{(2.3)}
$$
of two scalar differential operators
$$
L^{KDV} = \frac{\partial ^2}{\partial x^2} + U,
\eqno{(2.4)}
$$
$$
A =
\frac{\partial ^3}{\partial x^3} + \frac{3}{2}U \frac{\partial }
{\partial x} + \frac{3}{4} U_x.
$$
In such case we call that equation is represented by
$L,A$-pair.

The Novikov-Veselov equation as against of the KDV equation is represented
by $L,A,B$-triple
$$
\frac{\partial L}{\partial t} + [L,A] - BL = 0
\eqno{(2.5)}
$$
where
$$
L^{NV} = \partial \bar \partial  + U,
\eqno{(2.6)}
$$
$$
A = (\partial ^3 + V\partial) + (\bar \partial ^3 + \bar V \bar
\partial ) ,
$$
and
$$
B = \partial V+\bar \partial \bar V.
$$
Here  $B$ is the scalar operator of multiplication by function.

The representation of nonlinear equations
by $L,A,B$-triples was introduced by
Manakov  (\cite{M}) as a two-dimensional
generalization of representation
(2.3). Indeed, equation (2.3) preserves spectrum of
operator $L$ and deforms it's
eigenfunctions as follows :
$$
\frac{\partial \psi}{\partial t} = A \psi, \ \ L\psi=\lambda \psi.
\eqno{(2.7)}
$$
ŽOperator (2.6) is multidimensional and it's eigenspaces,
which correspond to
fixed eigenvalues, are generally
infinite-dimensional . Equation (2.5) deforms not all
eigenspaces but only
 the kernel of operator $L$ via the equation
$$
\frac{\partial \phi}{\partial t} = A \psi, \ \ L\phi=0.
\eqno{(2.8)}
$$

'In some sense the Novikov-Veselov equation is more natural two-dimensional
generalization of the KDV equation than the famous
Kadomtsev-Petviashvili equation
$$
(U_t - \frac{1}{4}(U_{xxx} + 6 U U_x))_x = \frac{3}{4} U_{yy}.
$$
This equation also reduces to the KDV equation in
the case when the function
$U$ does not depend on the space variable $y$ and it is represented by the
$L,A$-pair where the operator $L$ has the form
$$
\frac{\partial}{\partial y} - \frac{\partial ^2}{\partial x^2} -U(x,y).
$$
This operator differs from two--dimensional operator (2.6)  which is a usual
two-dimensional Schr\"odinger operator, i.e. the most natural two-dimensional
generalization of the one--dimensional Schr\"odinger operator (2.4).

‡ Notice that there exist two different deformations, of operator (2.6), which
have form (2.5):
$$
\frac{\partial L^{NV}}{\partial {t^\pm}} + [L^{NV},A^{\pm}] -
B^{\pm}L=0.
$$
These deformations are represented by $L,A,B$-triples for operators
$$
A^+ = \partial ^3 + V\partial ,
\bar \partial V = 3\partial U, B^+ = \partial V
$$
and
$$
A^- = \bar \partial ^3 + \bar V\bar \partial ,
\partial \bar V = 3\bar \partial U,
B^+ = \bar \partial \bar V.
$$
But these deformations do not preserve real potentials. In it's turn
equation (2.1), which, in fact,
is their linear superposition, transforms real potentials $U$ into
real ones.

\subhead
2.2. Modified Novikov-Veselov equation
\endsubhead

"There exists another famous integrable $1+1$-dimensional integrable equation
which is called modified Korteweg-de Vries equation (mKDV) :
$$
U_t = U_{xxx} + 24U^2 U_{x}.
\eqno{(2.9)}
$$
This equation is represented by  $L,A$-pair which we will not give here.
We only mention that $L$--operator has the following form
$$
L^{mKDV}=
\frac{\partial }{\partial x} -
\frac{1}{2}
\Big(\matrix
-1 &  4U \\
- 4U & 1
\endmatrix\Big).
\eqno{(2.10)}
$$

'Bogdanov introduced in  \cite{Bg1}  two-dimensional generalization
of the mKDV equation -- the modified Novikov-Veselov (mNV) equation
$$
U_t = (U_{zzz} + 3U_z V + \frac{3}{2}U V_z) +
(U_{{\bar z}{\bar z}{\bar z}} + 3U_{\bar z}{\bar V} +
\frac{3}{2}U{\bar V}_{\bar z})
\eqno{(2.11)}
$$
where
$$
V_{\bar z} = (U^2)_z.
\eqno{(2.12)}
$$

This equation is also, as the NV equation, a linear superposition of two
deformations of form (2.5) which are represented by $L,A,B$-triples
with common operator $L$ defined by
$$
L^{mNV}=
\Big(\matrix
\partial & -U \\
U & \bar \partial
\endmatrix\Big),
\eqno{(2.13)}
$$
 and the following $A$- and $B$-operators
$$
A^+ =
\partial^3 +
3
\Big(\matrix
0 & -U_z \\
0 & V
\endmatrix\Big)
\partial
+
\frac{3}{2}
\Big(\matrix
0 & 2UV \\
0 & V_z
\endmatrix\Big),
\eqno{(2.14)}
$$
$$
B^+=
3
\Big(\matrix
0 & U_z \\
-U_z & 0
\endmatrix\Big)
\partial
+
3
\Big(\matrix
0 & -UV \\
-U_{zz}-UV & 0
\endmatrix\Big)
$$
and
$$
A^- =
{\bar \partial}^3 +
3
\Big(\matrix
{\bar V}{\bar \partial} & 0 \\
U_{\bar z} & 0
\endmatrix\Big)
{\bar \partial}
+
\frac{3}{2}
\Big(\matrix
{\bar V}_{\bar z} & 0 \\
-2U{\bar V} & 0
\endmatrix\Big),
\eqno{(2.15)}
$$
$$
B^-=
3
\Big(\matrix
0 & U_{\bar z} \\
-U_{\bar z} & 0
\endmatrix\Big)
{\bar \partial}
+
3
\Big(\matrix
0 & U_{{\bar z} {\bar z}} + U{\bar V} \\
U{\bar V} & 0
\endmatrix\Big).
$$

These triples represent equations
$$
U_{t^+} = U_{zzz} + 3U_z V + \frac{3}{2}UV_z
$$
and
$$
U_{t^-} = U_{{\bar z}{\bar z}{\bar z}} + 3U_{\bar z}{\bar V} +
\frac{3}{2} U {\bar V}_{\bar z}
$$
where the function $V$ is defined by formula (2.12).

Analogously to the case of the NV equation, we can derive that

1) if the function $U$ depends only on one space variable $x$
when modified NV equations reduce to the mKDV equation ;

2) equation (2.11) transforms real potentials into real ones as against to
equations represented by $L^{mNV},A^{\pm},B^{\pm}$-triples ;

3) modified Novikov-Veselov equations deform the kernel of operator $L$
via the equations
$$
\frac{\partial \psi}{\partial t^{\pm}} = A^{\pm}\psi, \ \ \
L^{mNV} \psi = 0,
\eqno{(2.16)}
$$
and  deformation of eigenfunctions of $L^{mNV}$ via (2.11) is defined by
$$
\frac{\partial }{\partial t}
\Big(\matrix
\psi _1 \\
\psi _2
\endmatrix\Big)
=
\Big(
\partial^3 +
{\bar \partial}^3 +
3
\Big(\matrix
0 & -U_z \\
0 & V
\endmatrix\Big)
\partial
+
3
\Big(\matrix
{\bar V}{\bar \partial} & 0 \\
U_{\bar z} & 0
\endmatrix\Big)
{\bar \partial}
+
$$
$$
\frac{3}{2}
\Big(\matrix
0 & 2UV \\
0 & V_z
\endmatrix\Big)
 +
\frac{3}{2}
\Big(\matrix
{\bar V}_{\bar z} & 0 \\
-2U{\bar V} & 0
\endmatrix\Big)
\Big)
\Big(\matrix
\psi _1 \\
\psi _2
\endmatrix\Big).
\eqno{(2.17)}
$$

\subhead
2.3. Hierarchies of equations
\endsubhead

ŽOne of outstanding properties of equations
integrable by the inverse scattering
method is that they are included into
hierarchies of such equations which
are recursively defined.

 For instance, the  KDV equation and it's
modification ( mKDV)  are only
first members (for $k=1$)
of hierarchies of equations of the form
$$
U_{t_{2k+1}} = N_{2k+1}(U)
$$
where $N_{2k+1}(U)$ are nonlinear operators.
These equations are represented
by $L,A$--pairs with operators $L^{KDV}$ and $L^{mKDV}$ respectively.
For the KDV hierarchy operators $A$ have the following form
$$
A^{KDV}_k = \frac{\partial ^{2k+1}}{\partial x^{2k+1}} + ...
$$
where we denote by dots terms of lower orders.
These terms are defined by condition that
commutators of operators $L$ and $A$
would be operators of multiplication by scalars.

Thus we can say that the KDV hierarchy is attached to the operator
$L^{KDV}$. Analogously the mKDV hierarchy is defined.

The NV equation and it's modification are also included in hierarchies for
which $A$-operators take forms
$$
A^{NV}_k = \partial ^{2k+1} + ...
$$
and
$$
A^{mNV}_k =
\Big(\matrix
1 & 0 \\ 0 & 1
\endmatrix\Big)\partial ^{2k+1} + ... .
$$
respectively.

We also can say that these hierarchies are attached to operators
$L^{NV}$ (i.e., the two-dimensional Schr\"odinger operator) and
$L^{mNV}$ (i.e., the Dirac operator), respectively.

'In the soliton theory the method of defining of hierarchies by using of
so-called recursion operators is well-known.
For instance, the $k$-th equation of the KDV hierarchy takes the form
$$
U_{t_{2k+1}} = R^k(U_x)
$$
where the recursion operator is given by
$$
R = \frac{\partial ^2}{\partial x^2} + 3U + 3U_x (\frac{\partial }
{\partial x})^{-1}.
$$

ŠRegretfully, $2+1$-equations and the NV equations among them do not admit such
simple representation in terms of local operator $R$. There exists method
based on
using of bilocal operators (\cite{FS}) but it's realization is more
difficult than in the case of the KDV equations.

'However one can be confirmed in that by forms of higher equations.
We will describe, for instance, only second equations of these hierarchies.

The  NV2  equations
$$
U_{t^+_3} = \Phi_{NV2}(U), \ \
U_{t^-_3} = {\overline{\Phi_{NV2}(U)}}
$$
where
$$
\Phi(U)_{NV2} = \partial^5 U + V \partial^3 U + 2 V_z \partial^2 U +
(W + V_{zz})\partial U + W_z U =
\partial^5 U + \partial (V\partial^2 U + V_z \partial U + WU),
$$
are represented by $L,A,B$-triples with the following operators
$$
A^+= \partial^5 + V\partial^3 + V_z \partial^2 +  W \partial,
$$
$$
B^+= V_z \partial^2 + V_{zz} \partial + W_z,
$$
$$
A^- = \bar A^+, B^- = \bar B^+,
$$
and
$$
{\bar \partial } V = 5\partial U,
$$
$$
{\bar \partial }W = 5 \partial^3 U + 3 V\partial U + V_z U.
$$

The second equations of the mNV hierarchy are more complicated :
$$
U_{t^+_3} = \Phi_{mNV2}(U), \ \ \ \
U_{t^-_3}= {\overline{\Phi_{mNV2}(U)}},
$$
where
$$
\Phi(U)_{mNV2}=
U_{zzzzz} + 5 V U_{zzz} + \frac{15}{2}V_z U_{zz} +
5(V^2 -\frac{3}{2}V_{zz} + W) U_z +
5(VV_z - V_{zzz} + \frac{1}{2}W_z)) U
$$
and
$$
V_{\bar z} = (U^2)_z , W_{\bar z} = (U^2 V -U^2_z)_z.
$$

ŽOperators $A$ and $B$ are given by
$$
A^+ =
\partial ^5 +
5
\Big(\matrix
0 & -U_z \\
0 & V
\endmatrix\Big)
\partial ^3
+
5
\Big(\matrix
0 & UV - U_{zz} \\
0 & \frac{3}{2}V_z
\endmatrix\Big)
\partial ^2
+
$$
$$
5
\Big(\matrix
0 & \frac{1}{2}UV_z - U_z V - U_{zzz} \\
0 & V^2 - \frac{3}{2}V_{zz} + W
\endmatrix\Big)
\partial
+
5
\Big(\matrix
0 & UV^2 - 2 UV_{zz} + U_{zz}V + \frac{1}{2}U_z V_z + UW \\
0 & VV_z - V_{zzz} + W_z
\endmatrix\Big),
$$
$$
B^+=
5
\Big(\matrix
0 & U_z \\
-U_z & 0
\endmatrix\Big)
\partial ^3
+
5
\Big(\matrix
0 & U_{zz} - UV \\
-UV - 2U_{zz} & 0
\endmatrix\Big)
\partial ^2
+
$$
$$
5
\Big(\matrix
0 & U_z V + U_{zzz} - \frac{1}{2}UV_z \\
-\frac{3}{2}UV_z - 2U_{zzz} -3U_z V & 0
\endmatrix\Big)
\partial
-
$$
$$
5
\Big(\matrix
0 &  U(V^2 +W - 2V_{zz}) + U_{zz}V + \frac{1}{2}U_z V_z  \\
U(V^2+W-\frac{3}{2}V_{zz}) + U_{zz}V + 3 U_z V_z +
U_{zzzz} & 0
\endmatrix\Big),
$$
and
$$
A^-=
{\bar \partial} ^5 +
5
\Big(\matrix
{\bar V} & 0 \\
U_{\bar z} & 0
\endmatrix\Big)
{\bar \partial} ^3 +
5
\Big(\matrix
\frac{3}{2}{\bar V}_{\bar z} & 0 \\
-U {\bar V} + U_{{\bar z}{\bar z}} & 0
\endmatrix\Big)
{\bar \partial} ^2 +
$$
$$
5
\Big(\matrix
{\bar V}^2 - \frac{3}{2}{\bar V}_{{\bar z}{\bar z}} + {\bar W} & 0 \\
-\frac{1}{2}U{\bar V}_{\bar z}
+ U_{\bar z}V +
U_{{\bar z}{\bar z}{\bar z}} & 0
\endmatrix\Big)
{\bar \partial } +
5
\Big(\matrix
{\bar V}{\bar V}_{\bar z} - {\bar V}_{{\bar z}{\bar z}{\bar z}} +
{\bar W}_{\bar z} & 0 \\
-U{\bar V}^2 + 2U {\bar V}_{{\bar z}{\bar z}} -
U_{{\bar z}{\bar z}}{\bar V} - \frac{1}{2}
U_{\bar z}{\bar V}_{\bar z} -
U{\bar W} & 0
\endmatrix\Big),
$$
$$
B^-=
5
\Big(\matrix
0 & U_{\bar z} \\
-U_{\bar z} & 0
\endmatrix\Big)
{\bar \partial} ^3
+
5
\Big(\matrix
0 & U{\bar V}+ 2U_{{\bar z}{\bar z}} \\
-U_{{\bar z}{\bar z}} + U{\bar V} & 0
\endmatrix\Big)
\partial ^2
+
$$
$$
5
\Big(\matrix
0 &  \frac{3}{2}U{\bar V}_{\bar z} + 2U_{{\bar z}{\bar z}{\bar z}}
+ 3U_{\bar z}{\bar V} \\
-U_{\bar z}{\bar V} - U_{{\bar z}{\bar z}{\bar z}}
+ \frac{1}{2}U{\bar V}_{\bar z} & 0
\endmatrix\Big)
\partial
+
$$
$$
5
\Big(\matrix
0 & U({\bar V}^2 + {\bar W}-\frac{3}{2}{\bar V}_{{\bar z}{\bar z}})
+ U_{{\bar z}{\bar z}}{\bar V} + 3 U_{\bar z}{\bar V}_{\bar z} +
U_{{\bar z}{\bar z}{\bar z}{\bar z}} \\
U({\bar V}^2 +{\bar W} - 2{\bar V}_{{\bar z}{\bar z}}) +
U_{{\bar z}{\bar z}}{\bar V} + \frac{1}{2}U_{\bar z}{\bar V}_{\bar z}
& 0
\endmatrix\Big)
$$
"
Equations
$$
U_t= \Phi_{NV2}(U) + {\overline{\Phi_{NV2}(U)}}
$$
and
$$
U_t = \Phi_{mNV2}(U) + {\overline{\Phi_{mNV2}(U)}}
$$
preserve reality of potentials analogously to the case of the first equations.

\head
3. Weierstrass representation
\endhead

\subhead
3.1. Construction of minimal surfaces
\endsubhead

 
The most general method of constructing
minimal surfaces in the three-dimensional
Euclidean space was introduced by
Weierstrass, and we will start with it our
explanation of representation of surfaces.

Let take a pair of functions
$(\psi _1, \psi_2 )$ such that one of them,
$\psi_1$, is antiholomorphic and another, $\psi_2$, is holomorphic.
Let us suppose that these functions
are defined at the same simply connected
domain $S$ in a complex plane.
We have a system of equations
$$
\cases
\psi_{1z}=0, & \\
\psi_{2\bar z}=0. &
\endcases
\eqno{(3.1)}
$$

ŽLet us now define in terms of these functions a mapping
$$
T: S \rightarrow {\bold R}^3
\eqno{(3.2)}
$$
by the following formulas
$$
z \in S \rightarrow (X^1(z,\bar z),X^2(z,\bar z),X^3(z,\bar z)) \in
{\bold R}^3
$$
where
$$
X^1+iX^2=
i\int_{\gamma}({\bar\psi}^2_1 dz' - {\bar\psi}^2_2 d{\bar z}'),
$$
$$
X^1-iX^2=
i\int_{\gamma}(\psi^2_2 dz' - \psi^2_1 d{\bar z}'),
\eqno{(3.3)}
$$
$$
X^3=-\int_{\gamma}(\psi_2{\bar\psi}_1 dz' +
\psi_1{\bar \psi}_2 d{\bar z}').
$$

'Everywhere we suppose that integrals are taken over any path
$\gamma $ which lies in the domain $S$ and connects point $z$ with some initial
point $z_0$.
It follows from (3.1) that
integrands are closed forms and hence values of
integrals do not depend on choice of path $\gamma$.

Weierstrass had shown that

{\it surface $T(S)$ is minimal that means
that it's mean curvature vanishes
everywhere.}

\subhead
3.2. Generalized Weierstrass formulas
\endsubhead

…It is naturally to ask when formulas (3.3) define a surface in the
three-dimensional Euclidean space.
As one can see integrands ought to be closed forms and this condition
is sufficient. In the case of the Weierstrass representation that
follows from (3.1).

It turns out that if functions $\psi_1,\psi_2$
satisfy more general system
$$
\cases
\psi_{1z}=U\psi_2, & \\
\psi_{2\bar z}=-U\psi_1 &
\endcases
\eqno{(3.4)}
$$
with real potential $U$
then integrands in (3.3) occur to be closed forms.
Hence in this case formulas (3.3) define
a surface for every solution to system
(3.4).

That was shown in  \cite{Kon} where formulas
for induced metric and curvatures
were also derived.
Let us explain them here. Coordinates $(z,\bar z)$ are conformal and
in terms of them the metric tensor takes the form
$$
D(z,{\bar z})^2 dz d{\bar z}
\eqno{(3.5)}
$$
where
$$ D(z,\bar z) = |\psi_1 (z,\bar z)|^2 + |\psi_2 (z,\bar z)|^2.
$$£
The  Gaussian curvature is given by
$$
K=-\frac{1}{D^2}\Delta \log D,
\eqno{(3.6)}
$$
 and the mean curvature takes the form
$$
H=\frac{2U}{D}.
\eqno{(3.7)}
$$

Š This representation is not new . For instance, it is given in survey
\cite{Bb}, it was discussed by U. Abresch in middle 80's with
it's relation to constructing constant mean curvature surfaces,
in other terms it is given in a book of Eisenhart (\cite{E}),
and moreover it is equivalent to the well-known
Kenmotsu representation
Š(\cite{Ken}, see also \cite{HO}).

We will show in 3.4 that it is almost equivalent
to the definition of the second fundamental form (see Proposition 1).
"
Notice that the convenience of this form of representation is that
operator in linear problem (3.4) coincides with the operator $L^{mNV}$
to which the modified Novikov-Veselov hierarchy is attached.
That was the main source for definition of local deformation of surfaces which
was given in \cite{Kon} where this representation was rediscovered.

\subhead
3.3. On representation of surfaces by Weierstrass formulas
\endsubhead

 Let us consider the question how wide is the class of surfaces represented
by formulas (3.3 - 3.4).

Let
$$
F : \Sigma \rightarrow {\bold R}^3
\eqno{(3.8)}
$$
be a regular mapping of the domain $\Sigma $, of the complex plane  ${\bold C}$
with coordinates $(z,\bar z)$, into the three-dimensional Euclidean space,
and the induced metric is conformally Euclidean with respect to these
coordinates, i.e. a metric tensor takes the form
$D(z,\bar z)^2 dz d\bar z$.

'In this case the vector
$$
G(z) =
\Big(\frac{\partial F^1}{\partial z},
\frac{\partial F^2}{\partial z},\frac{\partial F^3}{\partial z}
\Big)
$$
satisfies evident equation
$$
\Big(\frac{\partial F^1}{\partial z}\Big)^2 +
\Big(\frac{\partial F^2}{\partial z}\Big)^2 +
\Big(\frac{\partial F^3}{\partial z}\Big)^2 = 0.
\eqno{(3.9)}
$$

That immediately follows from the following formula
$$
G(z) = \frac{\partial F}{\partial z} = \frac{1}{2}
\Big(\frac{\partial F}{\partial x}
-i\frac{\partial F}{\partial y}\Big)
\eqno{(3.10)}
$$
and the condition that the metric is conformally Euclidean :
$$
\Big(\frac{\partial F}{\partial x},\frac{\partial F}{\partial x}\Big)=
\Big(\frac{\partial F}{\partial y},\frac{\partial F}{\partial y}\Big),
\Big(\frac{\partial F}{\partial x},\frac{\partial F}
{\partial y}\Big)=0.
$$

The subvariety $Q \subset {\bold C}P^2$, which is defined in terms of
homogeneous coordinates
$(\varphi _1,\varphi _2,\varphi _3)$ by equation
$$
\varphi ^2_1 + \varphi ^2_2 + \varphi ^2_3 =0,
$$
is diffeomorphic to the Grassmann manifold
$G_{3,2}$ formed by two-dimensional subspaces of
${\bold R}^3$.
This diffeomorphism is given by the mapping
$$
G_{3,2} \rightarrow Q
$$
which corresponds to a plane, generated by a pair of orthogonal unit vectors
$(a_1,a_2,a_3)$ and  $(b_1,b_2,b_3)$, a point
$(a_1+ib_1,a_2+ib_2,a_3+ib_3) \in Q$.

Thus we can consider this mapping $G$ as the Gauss map.

 The Gauss map, defined in this manner, for surface (3.3) takes the form
$$
G(z) = (i({\bar \psi _1}^2 +\psi _2^2)/2,
({\bar \psi _1}^2 -\psi _2^2)/2, -\psi _2{\bar \psi _1}).
\eqno{(3.11)}
$$

'This formula gives us an idea how to prove the following Proposition.

\proclaim{Proposition 1}
‹Every regular conformally Euclidean immersion of surface into the
three-dimensional Euclidean space is locally defined by generalized
Weierstrass formulas (3.3 - 3.4).
\endproclaim

‡ For the sake of brevity, we did not mention that formulas
(3.3 - 3.4) represent locally every surface up to translation in
${\bold R}^3$. That is easy to see from (3.3).

"Proof of Proposition 1.

‡ We assume that $F^3_z \neq 0$ otherwise change coordinates in
${\bold R}^3$ to get that.

Let us compare '(3.10) and  (3.11) and define  functions
$$
\cases
\varphi _1 = \sqrt{F^2_{\bar z} + iF^1_{\bar z}},&\\
\varphi _2 = \sqrt{-(F^2_z+iF^1_z)}.&
\endcases
\eqno{(3.12)}
$$

It follows from ˆ(3.9) that
$$
F^3_z = -{\bar\varphi_1}\varphi _2.
$$

 Let us now remind the definition of the second fundamental form $h_{ij}$.
Let $D(z,\bar z)^2dz d{\bar z}$ be a metric tensor on surface (3.8).
We take in the tangent plane (at point $z$) an orthonormal basis
$$
e_1 = \frac{1}{D}\frac{\partial F}{\partial x}, \ \ \ \ \
e_2 = \frac{1}{D}\frac{\partial F}{\partial y}
$$
and extend it to a basis in ${\bold R}^3$ by adding a unit normal vector
$$
e_3 = e_1 \times e_2.
$$

ŠComponents of the curvature tensor are defined by the well-known decomposition
formulas (see, for instance, \cite{Ken}):
$$
\frac{\partial ^2 F}{\partial x^2} =
\frac{\partial D}{\partial x} e_1 -
\frac{\partial D}{\partial y} e_2 + D^2 h_{11} e_3,
$$
$$
\frac{\partial ^2 F}{\partial x \partial y} =
\frac{\partial D}{\partial y} e_1 +
\frac{\partial D}{\partial x} e_2 + D^2 h_{12} e_3,
$$
$$
\frac{\partial ^2 F}{\partial y2} =
-\frac{\partial D}{\partial x} e_1 +
\frac{\partial D}{\partial y} e_2 + D^2 h_{22} e_3.
$$

Substitute these expressions for second derivatives of $F$ into formulas for
$\varphi _{1z},\varphi _{2{\bar z}}$, derived from (3.12), and by
direct computations obtain
$$
\cases
\varphi_{1z} = \frac{DH}{2}\varphi _2,& \\
\varphi_{2\bar z} = -\frac{DH}{2}\varphi _1,&
\endcases
$$
where $H$  is a mean curvature.

Proposition 1 is proved.

'The important corollary of proposition 1 is the following Proposition.

\proclaim{Proposition 2}
Every regular analytic surface is represented by formulas
(3.3 - 3.4) globally.
\endproclaim

That follows from existence of local representation and unique analytic
continuation.

\subhead
3.4. Examples of surfaces
represented by Weierstrass formulas
\endsubhead

 Let us consider the simplest examples of surfaces represented by formulas
(3.3 - 3.4).

1) {\bf Surfaces of revolution.}

We assume, without loss of generality, that the axis  $OX^3$
is the axis of revolution.
In this case functions
$\psi _1$ and $\psi _2$ are given by
$$
\psi _1 = r_1(x) \exp{\frac{iy}{2}},\ \ \ \ \ \
\psi _2 = r_2(x) \exp{\frac{iy}{2}},
$$
 and system (3.4)  takes the form
$$
\Big(
\frac{\partial }{\partial x} -
\frac{1}{2}
\Big(\matrix
-1 & 4U \\
-4U & 1
\endmatrix\Big)
\Big)
\Big(\matrix
r_1 \\ r_2
\endmatrix\Big) = 0.
\eqno{(3.13)}
$$

‡Here a potential $U$ depends only on one variable $x$ and it is easy to
see that the matrix differential operator from linear problem (3.13)
coincides with the operator
$L^{mKDV}$ of form (2.10). Hence, in terms of the generalized Weierstrass
representation, the reduction of
'$L^{mNV}$ to $L^{mKDV}$ has natural geometrical meaning.

2) {\bf ‡ Closed surfaces with genus  $\geq 1$.}

Let  $F: \Sigma \rightarrow {\bold R}^3$ be
an immersion, of surface with genus
$g \geq 1$ , given by formulas  (3.3 - 3.4).

ˆIt is well-known that every closed oriented surface $\Sigma$
with positive genus is
uniformizable that means that there exists a mapping
$$
p: M \rightarrow \Sigma
$$
of simply connected surface $M$ with constant curvature
(the Euclidean plane for $g=1$ and the Lobachevsky plane for $g>1$)
which is conformal covering.

'In other words there exists a discrete subgroup
$\Gamma $ of a group of isometries of $M$ such that a factor--space
$M/\Gamma $
is conformally equivalent to the surface $\Sigma $.

 We consider cases $g=1$ and $g>1$ separately.

2.1) {\bf' Tori ($g=1$).}

'In this case a subgroup $\Gamma $ is
isomorphic to a free Abelian group with
rank 2 (i.e., two-dimensional lattice)
generated by a pair of independent shifts.

…If $\gamma \in \Gamma $ then
$$
\gamma ^*( dz d{\bar z}) = dz d{\bar z}
$$
and hence the following Proposition holds.

\proclaim{Proposition 3}
…Let $\Sigma $ be a two-dimensional torus immersed into
${\bold R}^3$ by formulas (3.3 - 3.4). Then
there exists  a lattice of periods
$\Gamma $, with rank 2, such that
potential $U(z)$ , metric tensor
$D(z)^2$, and mean curvature are invariant
with respect to action of $\Gamma$.
Functions $\psi _1,\psi _2$ at the same
time are transformed as follows
$$
\psi _1(z + \gamma ) = (\pm 1) \psi _1(z),
$$
$$
\psi _2(z + \gamma ) = (\pm 1) \psi _2(z),
$$
$$
z \rightarrow z + \gamma , \ \ \gamma  \in \Gamma.
$$
\endproclaim

2.2) {\bf Surfaces with genus $g>1$.}

'In this case a space $M$ is isometric to
the upperhalf plane
${\Cal H} = \{(x+iy) \in {\bold C} | y>0\}$
endowed with the  metric
$(dx^2+dy^2)/y^2$. ƒThe group of
isometries of ${\Cal H}$ is the group
$PSL(2,{\bold R})$ which acts by fractional
linear transformations
$$
z \rightarrow \frac{az+b}{cz+d},
$$
$$
a,b,c,d \in {\bold R}, \ \
ad-bc =1.
$$

"The action of elements of $PSL(2,{\bold R})$ on
differentials takes the form
$$
\gamma ^*(dz) = \frac{dz}{(cz+d)^2},\ \ \ \ \ \
\gamma =\Big(\matrix
a & b \\ c & d
\endmatrix\Big).
$$
ŽWe conclude that
$$
D(\gamma (z)) = |cz+d|^2 D(z).
\eqno{(3.14)}
$$
Since the mean curvature is invariant
($H(z)=H(\gamma (z))$), it follows from (3.7) that
$$
U(\gamma (z)) = |cz+d|^2 U(z).
\eqno{(3.15)}
$$

Now we are able to make the following conclusion.

\proclaim{Proposition 4}
Let a surface $\Sigma $, with genus $g>1$, is
immersed into
${\bold R}^3$ by formulas (3.3 - 3.4) and is
conformally equivalent to a surface
${\Cal H}/\Gamma $ where $\Gamma $ is a discrete subgroup of
$PSL(2,{\bold R})$.
Then 'metric tensor $D(z)^2$ and potential $U(z)$ are
transformed by elements of
$\Gamma$ by formulas  (3.14) and (3.15) respectively,
and $\Gamma$ acts
on functions $\psi _1$ and $\psi _2$ as follows
$$
\psi _1 (\gamma (z)) = (c{\bar z} + d) \psi _1(z),
$$
$$
\psi _2 (\gamma (z)) = (cz+d) \psi _2 (z).
$$
\endproclaim

\head
4. Deformation of surfaces  by the
modified Novikov-Veselov equation
\endhead

\subhead
4.1. Definition of deformation
\endsubhead

'In paper \cite{Kon} ŠKonopelchenko by using of
representation (3.3 - 3.4) defined
a new class of deformations of surfaces.
The mean observation of this paper is
that the operator from the linear problem (3.4)
coincides with the operator
$L^{mNV}$ to which the modified Novikov-Veselov
hierarchy is attached.
Hence the following deformation is naturally defined:

1) let  $F: S \rightarrow {\bold R}^3$ be
a surface immersed by formulas
(3.3 - 3.4) ;

2) assume that the potential $U(z,\bar z,t)$ is being transformed in $t$
via the modified Novikov-Veselov equation (2.11). At the same time
eigenfunctions  $\psi _1,\psi _2$ are being transformed via  equation (2.17)
and that generates deformation of surface $F$.

We call this deformation {\it modified Novikov-Veselov deformation} (mNV
deformation).

Moreover it is stated in \cite{Kon}
that every equation from the mNV hierarchy
generate deformation of this type.
Minimal surfaces correspond to zero potentials
and hence are stationary with respect to these flows.

An interesting observation was done in \cite{KT}
where the following Proposition was proved.

\proclaim{Proposition 5}
(\cite{KT})
1) An integral of squared mean curvature over
closed immersed (via (3.3 - 3.4))
surface $S$ , i.e. the value of
the Willmore functional at surface $S$,
is equal to
$$
W(S) = 4 \int_{F(S)} U(z,\bar z)^2 dz d{\bar z} \ \ \ ;
\eqno{(4.1)}
$$

2) If closed surface …$F(S)$ is deformed by
the mNV flow into closed surfaces and
lattices of periods of functions $U$ and $V$ are preserved then
the value of the Willmore functional
is also preserved.
\endproclaim

The proof of the first statement follows
immediately from formulas (3.5)
and (3.7). The second statement is
derived from the following
formula
$$
UU_t = (UU_{zz} - \frac{U_z^2}{2} + \frac{3}{2}U^2V)_z +
(UU_{{\bar z}{\bar z}} - \frac{U_{\bar z}^2}{2} + \frac{3}{2}
U^2{\bar V})_{\bar z},
\eqno{(4.2)}
$$
which itself follows from equation (2.11).

In the spirit of this Proposition it is natural to study global properties
of the mNV flow and it's relation to the theory of Willmore surfaces.
In \cite{KT} such investigation was started for tori of revolution.

We will not dwell here on results of \cite{KT} on tori of revolution and
explain here more general facts.

\subhead
4.2. Global deformations of closed surfaces
\endsubhead

In this subchapter we consider the question when the mNV flow transforms
closed surfaces into closed ones preserving their conformal structure.

First we thought that non-automorphic form of $L^{mNV}$ and $A$
operators (see Proposition 4) implies non-existence of
mNV-deformations , of surfaces with genus $g \geq 2$, that preserve
their closedness and conformal structure. But now
we are persuaded by F. Pedit and U. Pinkall that that strongly
depends on the correct understanding of constraint (2.12)
and definition of $V$. Probably at least deformation of periodic
Gauss map can be obtained in this manner.
This problem is still in question and thus we
will  discuss the case of surfaces with higher genus elsewhere.

Thus we restrict ourselves by deformations of tori.
First we prove the following Proposition.


\proclaim{Proposition 6}
'There exists procedure which uniquely corresponds to a smooth  double-periodic
potential $U(z)$ a function $V(z)$ which satisfies (2.12).
\endproclaim

Proof of Proposition 6.
"
Let $U(z)$ be a double-periodic function with a lattice of periods $\Gamma$
and
$\Gamma ^*$ be a lattice dual to $\Gamma$.

Any smooth function on the torus ‹${\bold C}/\Gamma = {\bold R}^2/\Gamma $
is decomposed into Fourier series with respect to basis formed by
eigenfunctions
of the operator $\bar \partial$. Notice that these functions also form a basis
of
eigenfunctions of the operator $\partial$.  These eigenfunctions are of the
form
$f(z|g^*) =\exp{2\pi i g^*(z)}$ where by $\gamma ^*(z)$
we mean a scalar product of $\gamma ^*$ and vector
$({\textstyle Re}\ z, {\textstyle Im}\ z)$.
For the sake of brevity we use only $z$ as argument but it is easy to notice
that these  functions are not holomorphic.

ŽIt is evident that for double-periodic function $w(z)$ there exists
double-periodic function $v(z)$ such that $v_{\bar z} = w$ if and only if
the Fourier series for $w(z)$
$$
w(z) = \sum_{g^* \in \Gamma^*} w_{\gamma^*} f(z|\gamma^*)
$$
does not contain terms corresponding to the
kernel of operator $\bar \partial$ (i.e., function $f(z|0)=1$).
In this case we can invert operator $\bar \partial$ evidently by using of
the Fourier decomposition.

If the function $w(z)$ is a derivative of the double-periodic
function itself then it's Fourier decomposition does not contain such term.
Let put $w(z)$ = $(U^2)_z$ and take a function $V(z)={\bar \partial}^{-1}w(z)$
uniquely determined by additional condition
$$
\int_{{\bold C}/\Gamma } V(z) dz d{\bar z} = 0.
$$
'This condition holds if and only if the Fourier series for $V(z)$ does not
contain terms which lie in the kernel of ${\bar \partial}$.

Proposition 6 is proved.

Notice that if we will  add to $V(z,t)$ a function which depends only on $t$
when
we will not change geometric deformation of surface but only include linear
translation of conformal coordinates $(z,\bar z)$.

Let us consider two integrals
$$
\frac{\partial(X^1(z,t)+iX^2(z,t))}{\partial t} =
2i \int \Omega _0
$$
and
$$
\frac{\partial X^3}{\partial t} =
-\int \Omega _1
$$
where
$$
\Omega _0 = \frac{1}{2}((\psi ^2_2)_t dz -(\psi ^2_1)_t d{\bar z}),
$$
$$
\Omega _1 = (\psi_{2t} {\bar \psi }_1 + \psi _2 {\bar \psi}_{1t}) dz +
(\psi_{1t} {\bar \psi }_2 + \psi_1 {\bar \psi}_{2t}) d{\bar z}.
$$

ŸExplicit formulas for differentials
$\Omega _0$ and $\Omega _1$ follows from (2.17). We omit them also as
rather large computations which need only formulas (2.12) and (3.3) and
give  the following result.

\proclaim{Proposition 7}
$$
1)\ \ \ \    \Omega _0 = d(f_1+g_1+f_2+g_2)
$$
where
$$
f_1 = \frac{3}{2}V\psi ^2_2 ,\ \ \ \ \
g_1 = \psi _2 \partial ^2\psi _2 - \frac{(\partial \psi _2)^2}{2},
$$
$$
f_2 = \frac{3}{2}{\bar V}\psi ^2_1 ,\ \ \ \ \
g_1 = \psi _1 {\bar \partial} ^2\psi _1 -
\frac{(\bar \partial \psi_1)^2}{2} ;
$$
$$
2) \ \ \ \    \Omega_1 = d(h_1+h_2)
$$
where
$$
h_1 =
{\bar \psi}_1\partial ^2\psi _2 +
\psi _2 \partial^2{\bar \psi }_1 -
\partial \psi _2 \partial {\bar \psi }_1 + 3V{\bar \psi }_1\psi _2,
$$
$$
h_2 =
\psi_1{\bar \partial}^2{\bar \psi}_2 +
{\bar \psi} _2 {\bar \partial}^2 \psi _1 -
{\bar \partial}{\bar \psi}_2 {\bar \partial}\psi_1 +
3{\bar V}\psi_1{\bar \psi}_2.
$$
\endproclaim

‡ 
Moreover two modified Novikov-Veselov deformations generated by $L,A,B$-triples
(2.14) and (2.15) satisfy to the following formal equations
$$
\frac{\partial(X^1(z,t^+)+iX^2(z,t^+))}{\partial t^+} =
2i \int d(f_1+g_1) ,
$$
$$
\frac{\partial(X^1(z,t^-)+iX^2(z,t^-))}{\partial t^-} =
2i \int d(f_2+g_2) .
$$
Žwhich can be useful for proving of analogues of Proposition 7 for
deformations generated by higher equations of the mNV hierarchy.

Now we can formulate the main Theorem.

'
\proclaim{Theorem 1}
…Let $\Sigma $ be a two-dimensional torus represented by formulas (3.3 - 3.4 )
with double-periodic potential $U(z)$, and let $U(z,t)$ be a solution to
equation (2.11) with
initial data $U(z,0) = U(z)$ and double-periodic potential $V(z,t)$. Then
the mNV flow deforms torus $\Sigma$ into tori $\Sigma_t$ which
are represented by (3.3 - 3.4) with potentials $U(z,t)$ , conformally
equivalent to $\Sigma$ and have the same value of the Willmore functional.
\endproclaim

Proof of Theorem 1.

By Proposition 7, forms $\Omega _0$ and  $\Omega _1$ are exact on torus
${\bold C}/\Gamma $ being differentials of double-periodic functions.
Therefore, a lattice of periods, which determines conformal class,
is preserved by the mNV flow.

Now it follows from Proposition 5 that a value of the Willmore functional is
also preserved.

Theorem 1 is proved.

In analytic case the stronger theorem holds.

\proclaim{Theorem 2}
TŒhe modified Novikov-Veselov equation induces via formulas (2.17) and
(3.3 - 3.4) deformation of immersed analytic tori. Moreover this deformation
preserves their conformal structures and values of the Willmore functional.
\endproclaim

"Proof of Theorem 2.

It follows from Proposition 2 that every analytic torus is
represented by formulas (3.3 - 3.4).
Since tori are analytic and by Proposition 6 and the Cauchy-Kowalewski theorem
a solution, of the modified Novikov-Veselov equation, which satisfies
conditions of Theorem 1 exists at least for small $t$.

Now Theorem 2 follows from Theorem 1.

Theorem 2 is proved.

\subhead
4.3. Clifford torus as stationary point of the
mNV flow
\endsubhead

ˆ

It follows from the definition of the potential $U(z)$ (see (3.7)) that
geometrically stationary points of the mNV flow, i.e., surfaces which are
transformed into images of itself with respect to translations in ${\bold
R}^3$,
correspond to stationary solutions of the mNV equation (2.11).

It is also naturally to expect that the simplest stationary solutions will be
one-dimensional, i.e. stationary solutions of the modified Korteweg--de Vries
equation.

We will show that the simplest stationary solution is realized by a prominent
surface, {\it Clifford torus}.

'Let $S^3$ be a unit sphere in the four-dimensional Euclidean space
${\bold R}^3$ with
coordinates $(x_1,x_2,x_3,x_4)$.
'The Clifford torus (in ${\bold R}^4$) is the image of the following
embedded torus
$$
{\bold R}^2 \rightarrow S^4 :
(x,y) \rightarrow (\frac{\cos{y}}{\sqrt{2}},
\frac{\sin{y}}{\sqrt{2}}, \frac{\cos{x}}{\sqrt{2}},
\frac{\sin{x}}{\sqrt{2}}).
\eqno{(4.3)}
$$

Let us consider the stereographic projection of
$S^4$ onto the plane  $x^4=-1$ with the pole $(0,0,0,1)$:
$$
(x_1,x_2,x_3,x_4) \rightarrow (\frac{-2x_1}{x_4-1},\frac{-2x_2}{x_4-1},
\frac{-2x_3}{x_4-1},-1).
\eqno{(4.4)}
$$

We call the image of the Clifford torus with respect to this projection
Clifford again.

The variables $(x,y)$ occur to be conformal and the metric tensor takes the
form
$$
\frac{4}{(\sqrt{2} - \sin{x})^2}( dx^2 + dy^2 ).
\eqno{(4.5)}
$$
ƒ Gaussian and mean curvatures are given by
$$
K= \frac{\sqrt{2}\sin{x} - 1}{4}, \ \
H=\frac{\sin{x}}{2\sqrt{2}}.
\eqno{(4.6)}
$$

Let us determine potential $U(x)$ by formula (3.7) :
$$
U(x) = \frac{\sin{x}}{2\sqrt{2}(\sqrt{2}-\sin{x})}.
\eqno{(4.7)}
$$

‹It follows from direct computations  that this potential induces the Clifford
torus by formulas (3.3 - 3.4). Let us also notice that potential (4.7)
satisfies
the following equation
$$
U_x^2 = -4U^4 + 2U^2 + \frac{U}{\sqrt{2}} + \frac{1}{16}.
\eqno{(4.8)}
$$

…If a solution of the mNV equation depends only on variable
$x - {\textstyle const}\cdot t$ then it satisfies to equation
$$
(U_{xxx} + 24 U^2 U_x - {\textstyle const}\cdot U_x)_x =0.
\eqno{(4.9)}
$$

ˆIf follows from (4.8) that the Clifford torus (4.7) satisfies (4.9).
Hence we conclude that the Clifford torus is a geometrically stationary
point of the mNV flow.

\head
5. Willmore functional
\endhead

'We already mentioned above (see Proposition 5 and Theorems 1 and 2) that the
mNV flow preserves values of the Willmore functional and briefly gave
definition of this functional.
Last years this functional attracted a lot of attention of geometers
(see history of it's investigation and explanation of a lot of facts about it
in
\cite{Wm}, also see \cite{ST,W,LY,Br,LS,Kus,FPPS,HJP,Sm,BBb,B}).

In this chapter we will try to get a brief survey of the modern history of the
Willmore conjecture and consider it's relation to the mNV flow.

Let $F: S \rightarrow {\bold R}^3$ be an immersed surface.
A value of the Willmore functional at this surface if defined by the following
formula :
$$
W(S) = \int _S H^2 d\mu.
\eqno{(5.1)}
$$
‡Here $d\mu$  the Liouville measure with respect to the induced metric on $S$.

This functional is conformally invariant,
i.e., any conformal transformation of
the three-dimensional Euclidean space
transforms any immersed surface into
another one with the same value of the Willmore functional.

We call surface {\it Willmore} if this
surface is a critical point of the
Willmore functional. The Euler-Lagrange
equation for this functional has the
form
$$
\Delta H + 2H(H^2 - K) =0,
\eqno{(5.2)}
$$
where $\Delta $ is the Laplace-Beltrami operator on surface
(\cite{Wm}).

ˆThe following Proposition can be obtained by using of direct computations.

\proclaim{Proposition 8}
…If surface is represented by formulas (3.3 - 3.4) then it is Willmore if and
only if the following equality holds
$$
\Delta U \cdot D - 2(U_xD_x+U_yD_y) + U \cdot \Delta D + 8U^3D=0
\eqno{(5.3)}
$$
where $z=x+iy$ and
$\Delta = \partial  ^2/\partial x^2 + \partial ^2/\partial y^2$.
\endproclaim

The simplest examples of Willmore surfaces are
stereographic projections of
minimal surfaces $M$ in $S^3$. Moreover an
area of a minimal surface
$M$ in $S^3$ is equal to
a value of the Willmore functional at it's
projection. All that was known to
Thomsen and Blaschke in 20-s. Conformal
properties of this functional and it's
relation to minimal surfaces gave Blaschke a
reason to call such surfaces
{\it conformally minimal} (\cite{Bl}).

But these examples do not cover the class of
Willmore surfaces (see,
for instance, \cite{P} where
the first examples of compact embedded
Willmore surfaces that are
not stereographic projections of minimal
surfaces in $S^3$ were obtained).

All Willmore spheres were classified by Bryant (\cite{Br})

ŽThe main attention attracts the conjecture posed by Willmore
in middle 60-s.

\proclaim{Willmore Conjecture}
"For immersed tori the Willmore functional satisfies the following inequality
$$
W \geq 2\pi^2,
\eqno{(5.4)}
$$
which is attained only for Clifford torus
and it's images under conformal transformations of ${\bold R}^3$.
\endproclaim

'
It's analogues for all genuses were posed in \cite{Kus} but this conjecture
is still open.

Simon proved that minimum is attained on an analytic minimal torus (\cite{Sm}).

ŠThe following list contains all known classes of tori for which Willmore
conjecture was proved.

1) 'In early 70-s Willmore and independently Shiohama and Takagi
(\cite{ST}) proved this conjecture for tube tori with constant radii.
Here we call torus tube if is  formed carrying a small circle round a closed
space curve such that the center moves along the curve and the plane of the
circle is the normal plane to the curve at each point.

2) •Hertrich-Jeromin and Pinkall (\cite{HJP})
generalized result of Willmore-Shio\-ha\-ma-Takagi
for tube tori with arbitrary radii,
i.e. radius of circle can vary along
the curve.

3) Langer and Singer (\cite{LS})
proved the Willmore conjecture for tori of
revolution.

4) "Li and Yau in paper \cite{LY}
that brings together the spectral
theory of the Laplace-Beltrami operator
with the theory of conformal invariants
proved this conjecture for tori whose
conformal structures are defined
by lattices generated by
vectors $(1,0)$ and $(a,b)$ where
$$
0 \leq a \leq \frac{1}{2}, \ \ \
\sqrt{1-a^2} \leq b \leq 1.
$$

'In terms of theta-functions all
Willmore tori are described in
\cite{BBb,B} (see also  \cite{FPPS}).
Regretfully theta-functional formulas are
very complicated and not rather
efficient for applications.

ŒWe propose the following conjecture.

\proclaim{Conjecture}
…Non-stationary, with respect to
the mNV flow, torus can not be a local minimum
of the Willmore functional.
\endproclaim

 As it seems to us this conjecture
looks truly because it is strange to expect
that minimum, of this variational
problem, taken up to conformal
transformations of ${\bold R}^3$ is degenerated.
Probably methods developed in \cite{W,Pm}
will be helpful for proving it.

If this conjecture is true then the
Willmore conjecture is reduced to
investigation of stationary points of the mNV flow. It is known from
the soliton theory that stationary solutions ought
to be simpler than general ones.
For instance, stationary solutions of
equations from the KDV hierarchy are
described by very simple hyperelliptic
functions. Of course the mNV equation is
$2+1$-equation and we can not expect for it so simple description.

'We also would like to pose the following question.

\proclaim{Question}
'Higher equations of the mNV hierarchy also  have first integrals.
What is a geometric meaning of critical points of these functionals ?
\endproclaim

'Similarity of formulas for mNV and mNV-2 equations shows that it needs to
expect that these flows will deform tori into tori.
Thus these deformations ought to have geometric meaning.
Most probably these flows preserve conformal structures and these flows have
origin in conformal geometry.

\enddocument

\enddocument